\documentclass[11pt]{article}
\usepackage{amssymb}
 \usepackage{amsthm}
\usepackage{float,graphicx,amsmath,amsfonts,amssymb,amsthm,color,slashbox,amsbsy,ragged2e}
\usepackage{hyperref}
 \usepackage{lineno}

\usepackage[figuresright]{rotating}

\def\Rea{\Bbb R}
\newcommand{\rb}{\boldsymbol{r}}

\newcommand{\teb}{\boldsymbol{\theta}}

\newcommand{\si}{\sigma}

\newcommand{\bq}{\begin{equation}}
\newcommand{\eq}{\end{equation}}
\newcommand{\bqs}{\begin{equation*}}
\newcommand{\eqs}{\end{equation*}}
\newcommand{\bqa}{\begin{eqnarray}}
\newcommand{\eqa}{\end{eqnarray}}
\newcommand{\bqas}{\begin{eqnarray*}}
\newcommand{\eqas}{\end{eqnarray*}}
\newcommand{\bc}{\begin{cases}}
\newcommand{\ec}{\end{cases}}
\newcommand{\bt}{\begin{thm}}
\newcommand{\et}{\end{thm}}

\title{A new simple and powerful normality test for   progressively Type-II censored data}
\author{ Hamzeh Torabi,  Sayyed Mahmoud  Mirjalili, Hossein Nadeb\\
Department  of Statistics, Yazd University,  Yazd, Iran,\\
}
\date{}

\begin{document}
\maketitle

\begin{abstract}
In this paper, a new goodness-of-fit test for a location-scale family  based on progressively Type-II censored order statistics 
is proposed.  
Using Monte Carlo simulation studies, the present researchers have observed that the proposed test for normality is consistent and quite powerful in
comparison with  existing goodness-of-fit tests based on progressively Type-II censored data.  Also, the new test statistic for a real data set is used and the results show that our new test statistic performs well. 
\end{abstract}


\section{Introduction}\label{}

One of the most interesting problems in statistics is finding a distribution which fits to a given set of data. In other words, it is desired to test whether a specific distribution coincides with given data or not. To review the classical goodness-of-fit test problem, let $X_1,...,X_n$ be random sample from an absolutely  continuous population  with cumulative distribution function (CDF) $F(.)$, and probability density function (PDF) $f(.)$. Based on the observed  sample $x_1,...,x_n$,  hypotheses testing of interest is
\begin{equation}\label{e1}
\left\{\begin{array}{l}H_0:~f=f_0\\ H_1:~f\neq f_0,\end{array}\right.
\end{equation}
where $f_0(x)=f_0(x;\teb)$, where $\teb\in \Theta\subseteq \Bbb R^k$ is a $k$-vector parameter for some  $k\in \Bbb N$. 
For more overview on the  topic of goodness-of-fit test,  refer to the books by D'Agostino and Stephens  \cite{agos} and Huber-Carol et al.  \cite{hub}.

Most of goodness-of-fit tests are based on the distance between empirical distribution function (EDF)  and theoretical distribution functions over the interval $(0, 1)$,  the null hypothesis is rejected if the distance is too large in some metrics. However,  one can construct a goodness-of-fit test  based on order statistics  in terms of the deviation of each order statistic $U_{i:n}$ from its expected value $i/(n+1)$, say, $V_i=U_{i:n}-i/(n+1)$. Statistics that can be considered in this regard are as follows: 
$$C^+_n=\mathop {\max }\limits_{1 \le i \le n}(V_i),  	C^-_n=\mathop {\max }\limits_{1 \le i \le n}(-V_i), C_n=\mathop {\max } (C^-_n,C^+_n), $$
$$K_n=C^-_n+C^+_n,   T^{(1)}=\sum\limits_{i = 1}^n {\frac{{V_i ^2 }}{n}}, T^{(2)}=\sum\limits_{i = 1}^n {\frac{{|V_i | }}{n}} $$
For the suitability of uniformity, the upper tail of the appropriate null distribution is usually used to test.  One may refer to Brunk  \cite{brunk}, Stephens  \cite{stephen} and Hegazy and Green  \cite{hegazy}, for more discussion on these statistics.

Goodness-of-fit testing can also be done based on the spacings $D_i=U_{i:n}-U_{i-1:n}$, $i=1,...,n+1$, where $U_{0:n}=0$ and $U_{n+1:n}=1$. However, several
statistics based on spacings  have been reported in the literature, including Greenwood's  \cite{greenwood} statistic $G_{(n)}=\sum\limits_{i = 1}^{n+1} D_i^2$, Quesenberry and Miller's  \cite{quesenberry} statistic  $Q=\sum\limits_{i = 1}^{n+1} D_i^2 + \sum\limits_{i = 1}^nD_iD_{i+1}$ and Moran's  \cite{moran} statistic $M(n)=-2\sum\limits_{i = 1}^{n+1}\ln\big((n+1)D_i\big)$.  The null hypothesis will be rejected for large values of these statistics. Also Torabi  \cite{to1,to2} has introduced a new and  general method for estimation and  hypotheses testing using spacing. 

The classical goodness-of-fit tests for complete data can no longer be used for progressively Type-II censored data, Pakyari and Balakrishnan  \cite{pak2013} employed a modification to the aforementioned statistics based on order statistics and spacings including the class of $C$ statistics, Greenwood's statistic, and Quesenberry and Miller's statistic, making them suitable for
progressively Type-II censored data. 

For progressive Type-II censoring, we refer to the recent survey paper by Balakrishnan   \cite{bala2007} and the monograph by Balakrishnan and Aggarwala   \cite{bala2000}. In progressive Type-II censoring, it is assumed that the removals of still operating units are
carried out at observed failure times and that the censoring scheme $({r_1 ,r_2 ,...,r_m })$
is known in advance. Moreover, the number of units ($n$) and the number of observed
failure times ($m$) are prefixed. Starting all $n$ units at the same time, the
first progressive censoring step takes place at the observation of the first failure time $X_{1:m:n}$, 
at this time, $ r_1$ units are randomly chosen from the still operating units and
withdrawn from the experiment. Then, the experiment continues with the reduced
sample size $n-r_1-1$. After observing the next failure at time $X_{2:m:n}$,  $r_2$ units are randomly removed from $n-r_1-2$ active units. This process continued until the $m$th failure is observed. Then, the experiment ends. The failure
times $ X_{1:m:n}, . . . , X_{m:m:n}$ are called progressively Type-II censored order statistics and $ x_{1:m:n}, . . . ,x_{m:m:n}$ are the corresponding observations. 
 For their relation to order statistics and other related models of order random variables, one may  refer to  Balakrishnan  \cite{bala2007}.

Goodness-of-fit test for the exponential distribution based on spacings from progressively Type-II censored data introduced by Balakrishnan et al.  \cite{bala2002}, then they   
extended their method to general location--scale families of distributions \cite{bala2004}. Also Wang  \cite{wang}
proposed another goodness-of-fit test for the exponential distribution under progressively Type-II
censored samples. Recently, Pakyari and Balakrishnan  \cite{pak2012} proposed a modification to the EDF
goodness-of-fit statistics under progressively Type-II censored data. One  may also
refer to  \cite{aho,balarad,bar,lin,lurie,pettit,smith,stephen} for some other developments in this regard. 

 In Section 2,  we review the test statistics based on spacings that are  modification to the previously defined $ C^+, C^-, C, K, T^{(1)}, T^{(2)}$, modification to the Greenwood's statistic and modification to the Quesenberry and Miller's statistic for the progressively Type-II censored data that proposed by Pakyari and Balakrishnan  \cite{pak2013}. In Section 3,  we propose a new test statistic that  will be used for test of normality under the progressively Type-II censored data. In Section 4, we investigate consistency of our test statistic using a simulation study under five progressively Type-II censored schemes. The power of the proposed test is then assessed through Monte Carlo simulations in Section 5, and its performance is compared with those of the test procedures introduced earlier by Balakrishnan et al.  \cite{bala2004} and Pakyari and Balakrishnan \cite{pak2013}. It is shown  that the proposed goodness-of-fit test to be more powerful than or at least as good as the tests of Balakrishnan et al.  \cite{bala2004}  and Pakyari and Balakrishnan  \cite{pak2013} for different choices of sample sizes and progressive censoring schemes. In Section 6, we illustrate the application of  proposed goodness-of-fit procedure with a real data set.
\section{ Review on the test statistics based on spacings}
In  Section 1, some several test statistics based on the deviation between order statistics and the corresponding expected value in the case of a complete sample were presented. These statistics were extended to progressively Type-II censored data by Pakyari and Balakrishnan \cite{pak2013} as follows
 $$C^+_{m:n}=\mathop {\max }\limits_{1 \le i \le m}(V_{i:m:n}),  	C^-_{m:n}=\mathop {\max }\limits_{1 \le i \le m}(-V_{i:m:n}), C_{m:n}=\mathop {\max } (C^-_{m:n},C^+_{m:n}), $$
$$K_{m:n}=C^-_{m:n}+C^+_{m;n},   T^{(1)}_{m:n}=\sum\limits_{i = 1}^m {\frac{{V_{i:m:n} ^2 }}{m}}, T^{(2)}_{m:n}=\sum\limits_{i = 1}^m {\frac{{|V_{i:m:n} | }}{m}}, $$ 
where in this case, $V_{i:m:n}=U_{i:m:n}-\mu_{i:m:n}$, which $U_{i:m:n}$ is the $i$th order statistic from uniform (0,1) distribution base on Type-II Progressive censored data and $\mu_{i:m:n}$ is its expected value, i.e 
$$\mu_{i:m:n}=1-\prod_{k=m-i+1}^m \frac{k+\sum_{j=m-i+1}^m r_j}{1+k+\sum_{j=m-i+1}^m r_j},\quad i=1,\ldots,m.$$ 
It is easy to show that all the above statistics are location--scale invariant. If the null hypothesis is true,
we expect that $V_{i:m:n}$ to be small and consequently the above test statistics to be small.
If the above test statistics exceed the corresponding upper-tail null critical values, the null hypothesis  may be rejected. 
Recently, several goodness-of-fit statistics based on spacings have been developed. The 
one-step spacings are defined by
$$S_i=(n-r_1-r_2-...-r_{i-1}-i+1)(U_{i:m:n}-U_{i-1:m:n}),\quad   i=1,2,...m,$$
where $U_{0:m:n}=0$. It was shown by Balakrishnan and Aggarwala  \cite{bala2000} and Viveros and Balakrishnan  \cite{bala2000} that if the underlying distribution is exponential, then $S_1,S_2,...,S_m$ are independent and identically distributed as exponential with the scale parameter $\sigma$.

The following statistics are based on the spacings that generalized by  Pakyari and Balakrishnan  \cite{pak2013} under the progressively Type-II censored data:
\begin{itemize}
\item 
Statistics based on the sum of squares of the spacings, which are the generalization of Greenwood's statistic for progressively Type-II censored samples, are simply of the form 
$$G_{m:n}=\sum\limits_{i = 1}^m S_i^2.$$
\item
 The generalization of Quesenberry and Miller's statistic for progressively Type-II
censored samples will be of the form
$$Q_{m:n}=\sum\limits_{i = 1}^m S_i^2 + \sum\limits_{i = 1}^{m-1} S_iS_{i+1}.$$
	 The exact distributions of $G_{m:n}$
and $Q_{m:n}$ are not available explicitly but by Monte Carlo simulations the  percentage points will be determined.

\item 
The above statistics can also be defined in terms of higher order spacings. The overlapping
$k$-step spacings, for integer $k$, are defined as
$$S_i^{(k)}=(n-r_1-r_2-...-r_{i-1}-i+1)(U_{i+k-1:m:n}-U_{i-1:m:n}),   i=1,2,...m,$$
	with $U_{l:m:n}$ for $ l > m$. Hartley and Pfaffenberger   \cite{hartly} presented that the higher order spacings could be useful for testing large complete
samples. The extensions of Greenwood's statistic and Quesenberry and Miller's statistic in terms
of overlapping k-spacings take the forms
$$G_{m:n}^{(k)}=\sum\limits_{i = 1}^m (S_i^{(k)})^2.$$

The null hypothesis of uniformity is rejected if these statistics are too large.
\item
Balakrishnan et al's\cite{bala2004} test statistic was defined as below: 
$$T=\frac{\sum\limits_{i = 2}^{m-1} (m-i)G_i}{(m-2)\sum\limits_{i = 2}^{m} G_i}$$
where 
$$G_i=\frac{S_i}{E(s_i)}=\frac{U_{i:m:n}-U_{i-1:m:n}}{\mu_{i:m:n}-\mu_{i-1:m:n}}$$
In the next section, we propose a new test statistic and in Section 5 compare it with the test statistics reviewed in this section. 
\end{itemize}
\section{Proposed test}
\begin{figure}[ht]
\centerline{\includegraphics[width=7cm]{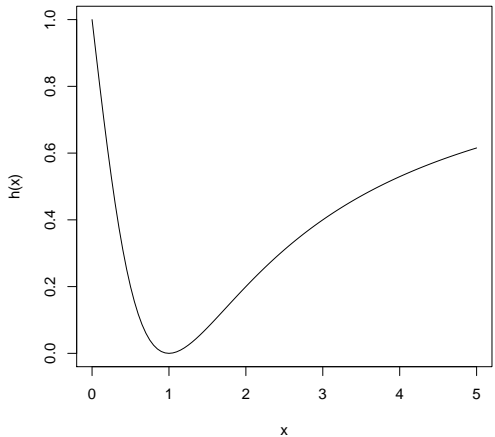}}
\caption{\small{ Plot of the function $h(x)$.}}
\label{fig}
\end{figure}

In this section, we propose a new  approach for goodness-of-fit testing for normality under progressively Type-II censored data. 
Consider again the goodness-of-fit testing problem (\ref{e1}) 
based on  $ X_{1:m:n}, . . . , X_{m:m:n}$, where  
  $f_0(x;\mu,\si)=1/\sqrt{2\pi \sigma^2} e^{-(x-\mu)^2/2\sigma^2}$, $x\in\mathbb{R}$, in which $\mu\in\mathbb{R}$ and $\sigma>0$ are both unknown. 
 Suppose $\hat \mu$ and $\hat \si$ are the MLEs of $\mu$ and $\si$ based on $ X_{1:m:n}, . . . , X_{m:m:n}$. Because of consistency of the ML estimators,   we expect $F_0(X_{i:m:n},\hat \mu, \hat \sigma)$ has the same distribution as  $U_{i:m:n}$; so it is justifiable  that $\frac{{F_0(X_{i:m:n},\hat \mu, \hat \sigma)}}{\mu_{i:m:n}} \simeq 1$.  Our proposed test is based on this ratio. More precisely, define
$${\rm H}_{m:n}=\frac1m\sum\limits_{i= 1}^m h\left(\frac{{F_0(X_{i:m:n},\hat \mu, \hat \sigma)}}{\mu_{i:m:n}}\right), $$
where   $h: (0, \infty) \rightarrow \Bbb R^+$ is  assumed to be continuous, decreasing
on $(0, 1)$ and increasing on $(1, \infty)$ with the absolute minimum at $x=1$ such that $h(1)=0$. 
In the simulation study for comparison of powers, by trying  some different  choices of $h$, the best choice is  
$$  h(x)=\frac{(x-1)^2}{x^2+1},$$ that has the maximal power.
Plot of the function $h$ is given in $Figure$ \ref{fig}.

We know that MLE of $\mu$ and $\sigma$ are location-scale invariant for $\mu$ and $\sigma$, respectively. Therefore under a location-scale transformation, the distribution of ${\rm H}_{m:n}$ does not depend on the parameters $\mu$ and $\si$ under location-scale transformations.

It is expected that  the null hypothesis of normality is rejected if the statistic  ${\rm H}_{m:n}$ is  too large; Thus the critical region is of the form ${\rm H}_{m:n}>c$, for some $c>0$. But for finding $c$ for a test of size $\alpha$, the exact distribution of ${\rm H}_{m:n}$ could not be  explicitly obtained, fortunately using Monte Carlo simulations, the  critical  points can be determined.
In Section 4, the consistency  of our proposed test  has checked and in Section 5 critical values of our test statistic has gained and its power has compared with the power of existence  test statistics  based on Monte Carlo simulations against Student's t, Logistic and Double Exponential models. 

An adequate  test statistic  for a goodness-of-fit test problem should be consistent, i.e, with increasing sample size, it is expected that the statistic tends to a finite value, especially under $H_0$ tends to zero. We cannot prove consistency of our test statistic but  using a Monte Carlo simulation study, is proved and discussed  in Section \ref{sec-con}. 
\section{Consistency of the new statistic using a simulation study}\label{sec-con}
In this section, we investigate consistency of our test statistic using a simulation study under five progressively Type-II
censored schemes. To illustrate the goal, we consider 5 various  censoring schemes as follows:
\begin{itemize}
\item {\bf Scheme 1}: a progressive Type-II censoring scheme with constant removal,   $\rb=(1,1,...,1)$, in this case $n=2m$;
\item {\bf Scheme 2}:  a progressive Type-II censoring scheme with  increasing  removal,   $r_i=i $ for $i=1,2,...,m$, in this case $n=m(m+3)/2$;
\item {\bf Scheme 3}: a progressive Type-II censoring scheme with decreasing  removal,   $r_i=m-i+1 $ for $i=1,2,...,m$, thus $n=m(m+3)/2$;
\item {\bf Scheme 4}: 
a Type-II censoring,   $r_i=0 $ for $i=1,2,...,m-1$, $r_m={m}/{5}$, hence $n=1.2m$;
\item {\bf Scheme 5}: complete data, i.e.,    $r_i=0 $ for $i=1,2,...,m$, thus $n=m$;
\end{itemize} 

As it is stated,  the normal model is considered as the parent model in $H_0$ but it can be changed with any location-scale model because of the structure of test statistic. Against this model, we consider some alternative models   as follows:
\begin{itemize}
\item Student's t distribution with $\nu$ degrees of freedom ($t_{(\nu)}$) with the density function
$$f(x;\nu)=\frac{\Gamma((\nu+1)/2)}{\sqrt{\nu\pi}\Gamma(\nu/2)}(1+\frac{x^2}{\nu})^{-(\nu+1)/2},\quad  \nu>0.$$
\item Logistic distribution with parameters $\mu$ and $\si$ ($L(\mu,\si)$) with the density function
$$f(x;\mu,\sigma)= \frac{\frac{1}{\sigma}\exp\left[\frac{-(x-\mu)}{\sigma}\right]}{\left(1-\frac{1}{\sigma}\exp\left[\frac{-(x-\mu)}{\sigma}\right]\right)^2},\quad x\in \Rea,\mu\in\mathbb{R},\sigma>0.$$
\item 
Double exponential distribution with parameters $\mu$ and $\si$ ($DE(\mu,\si)$) with the density function
$$f(x;\mu,\sigma)=\frac{1}{2\sigma}\exp\left[\frac{-|x-\mu|}{\sigma}\right],\quad x\in \Rea,\mu\in\mathbb{R},\sigma>0. $$
\end{itemize}
For more details of these distributions refer to Casella and Berger  \cite{Cas-Ber}.

Results  that are given in Tables 1-5,  show that under the standard normal distribution the values of our new statistic tend to zero when $m$ increases (and hence $m$ increases), but under the alternative distributions such as $t_{(3)}$, $t_{(4)}$, $L(0,1)$ and $DE(0,1)$  the values of our new statistic does tend to a non zero value for all five schemes.
\begin{table}[ht]
\begin{center}
\caption{Consistency of ${\rm H}_{m:n}$ using the Monte Carlo simulations for Scheme 1.}\label{t1}
\begin{tabular}{ccccccc}
\hline
$n$ &  $ m$  & $N(0,1)$  & $t_{(3)}$ &$t_{(4)}$ &$DE(0,1)$& $L(0,1)$\\
\hline
50& 25 & 0.0301 & 0.0617&0.0518&0.0580&0.0405 \\
100&50&0.0184&0.0541&0.0426&0.0470&0.0292\\
200&100&0.0112&0.0501&0.0366&0.0406 &0.0223\\
300&150&0.0084&0.0481&0.0345&0.0384 &0.0195\\
400&200&0.0067&0.0475&0.0334&0.0375&0.0180\\
500&250&0.0057&0.0459&0.0325&0.0368&0.0169\\
600&300&0.0050&0.0397&0.0314&0.0364&0.0134\\
                      \hline
\end{tabular} 
\end{center}
\end{table}
\begin{table}[ht]
\begin{center}
\caption{ Consistency of ${\rm H}_{m:n}$ using the Monte Carlo simulations for Scheme 2.}\label{t2}
\begin{tabular}{ccccccc}
\hline
$n$ &  $ m$  & $N(0,1)$  & $t_{(3)}$ &$t_{(4)}$ &$DE(0,1)$& $L(0,1)$\\
\hline
65&10&0.0639&0.0951&0.0865&0.0789&0.07356\\
230&20&0.0391&0.0818&0.0695&0.0535&0.0504\\
430&40&0.0241&0.0737&0.0594&0.0365&0.0346\\
1890&60&0.0178&0.0682&0.0546&0.0293&0.0280\\
3320&80&0.0144&0.0630&0.0525&0.0246&0.0244\\
5150&100&0.0123&0.0471&0.0464&0.0224&0.0217\\
                      \hline
\end{tabular} 
\end{center}
\end{table}
\begin{table}[ht]
\begin{center}
\caption{Consistency of ${\rm H}_{m:n}$ using the Monte Carlo simulations for Scheme 3.}\label{t3}
\begin{tabular}{ccccccc}
\hline
$n$ &  $ m$  & $N(0,1)$  & $t_{(3)}$ &$t_{(4)}$ &$DE(0,1)$& $L(0,1)$\\
\hline
65&10&0.0665&0.1124&0.0994&0.0949&0.0812\\
230&20&0.0449&0.1128&0.0924&0.0785&0.0636\\
430&40&0.0295&0.1191&0.0922&0.0647&0.0522\\
1890&60&0.0230&0.1251&0.0957&0.0579&0.0452\\
3320&80&0.0188&0.1249&0.0974&0.0540&0.0425\\
5150&100&0.0160&0.1195&0.0978&0.0503&0.0408\\
                      \hline
\end{tabular} 
\end{center}
\end{table}

\begin{table}[h]
\begin{center}
\caption{ Consistency of ${\rm H}_{m:n}$ using the Monte Carlo simulations for Scheme 4.}\label{t4}
\begin{tabular}{ccccccc}
\hline
$n$ &  $ m$  & $N(0,1)$  & $t_{(3)}$ &$t_{(4)}$ &$DE(0,1)$& $L(0,1)$\\
\hline
60&50&0.0170&0.0429&0.0343&0.0397&0.0247\\
120&100&0.0104&0.0391&0.0291&0.0343&0.0185\\
180&150&0.0077&0.0379&0.0269&0.0323&0.0158\\
240&200& 0.0064&0.0369 &0.0258 &0.0318&0.0145 \\
300&250& 0.0051&0.0361 &0.0247 &0.0312&0.0135 \\
360&300& 0.0046&0.0350 &0.0243 &0.0305&0.0130 \\
420&350& 0.0041&0.0313 &0.0239 &0.0300&0.0101 \\
480&400& 0.0037 &0.0215 &0.0193 &0.0274&0.0097 \\
                      \hline
\end{tabular} 
\end{center}
\end{table}
\begin{table}[ht]
\begin{center}
\caption{ Consistency of ${\rm H}_{m:n}$ using the Monte Carlo simulations for Scheme 5.}\label{t5}
\begin{tabular}{ccccccc}
\hline
$n$ &  $ m$  & $N(0,1)$  & $t_{(3)}$ &$t_{(4)}$ &$DE(0,1)$& $L(0,1)$\\
\hline
50&50&0.0165&0.0371&0.0297&0.0467&0.0215\\
100&100&0.0103&0.0351&0.0250&0.0265&0.0159 \\
200&200&0.0063&0.0334&0.0225& 0.0235& 0.0121\\
400&400&0.0037&0.0335&0.0206&0.0218 &0.0095\\
800&800&0.0022&0.0334&0.0196& 0.0208&0.0083\\
1600&1600&0.0013&0.0334&0.0189& 0.0203&0.0075\\
3200&3200&0.0007&0.0334&0.0189& 0.0203&0.0075\\
                      \hline
\end{tabular} 
\end{center}
\end{table}
\section{Simulation Study}
In this section,  we assess the power of the our new statistic by comparing the simulated power values
with those of the test of Balakrishnan et al. \cite{bala2004} and  Pakyari and Balakrishnan \cite{pak2013}. We calculated the power of the proposed test for testing of normality against some different alternatives with simulating  10,000 random samples for some different choices of sample sizes and progressive
censoring schemes. For comparative purposes, all 27 censoring
schemes used by Balakrishnan et al.   \cite{bala2004} and  Pakyari and  Balakrishnan   \cite{pak2013} in their studies  are  considered again here, and these are listed in Table \ref{t6}. Also the simulated critical values of ${\rm H}_{m:n}$ for every 27 censoring scheme has listed in Table \ref{t10}. All the simulations were carried out in $\tt R$ software.

\begin{table}
\begin{center}
\caption{ Progressive censoring schemes used in the Monte Carlo simulations. }\label{t6}
\begin{tabular}{cccc}
\hline
Scheme no. &  $ n$  &$m$  & $\rb=(r_1,r_2,...,r_m)$\\
\hline
$[1]  $      &20 & 8 & $ r_1=12, r_i=0 \mbox{ for } i\ne1    $\\
  $[2] $       &20 & 8 & $ r_8=12, r_i=0 $ for $i\ne8$\\
  $[3] $           &20 & 8 & $ r_1=r_8=6, r_i=0$ for $i\ne1,8    $\\
&&&\\
$[4]$        &20 & 12 & $r_1=8, r_i=0 $ for $i\ne8$\\
$[5]   $         &20 & 12 & $ r_{12}=8, r_i=0$ for $ i\ne12    $\\
$[6]$        &20 & 12 & $r_3=r_5=r_7=r_9=2, r_i=0$ for $ i\ne3,5,7,9$\\
&&&\\
$[7]   $         &20 & 16 & $ r_1=4, r_i=0$ for $ i\ne1    $\\
$[8]$        &20 & 16 & $r_{16}=4, r_i=0  $ for $ i\ne16 $\\
$[9]   $         &20 & 16 & $ r_5=4, r_i=0$ for $i\ne5    $\\
&&&\\
$[10]$        &40 & 10 & $ r_1=30, r_i=0$ for $i\ne1$\\
$[11]   $         &40 & 10 & $ r_{10}=30, r_i=0$ for $ i\ne10    $\\
$[12]$        &40 & 10 & $r_1=r_5=r_{10}=10, r_i=0$ for $ i\ne1,5,10$\\
&&&\\
$[13] $           &40 & 20 & $ r_1=20, r_i=0$ for $i\ne1   $\\
$[14]$        &40 & 20 & $r_{20}=20, r_i=0 $ for $ i\ne20$\\
$[15]  $          &40 & 20 & $ r_i=1$, for $ r_i=1,2,...,20    $\\
&&&\\
$[16]$        &40 & 30 &$ r_1=10, r_i=0 $ for $ i\ne1$\\
$[17]   $         &40 & 30 & $ r_{30}=10, r_i=0$ for $i\ne30   $\\
$[18]$        &40 & 30 & $r_1=r_{30}=5, r_i=0$ for $ i\ne1,30$\\
&&&\\
$[19]$        &60 & 20 & $r_1=40, r_i=0  $ for $ i\ne1$\\
$[20]   $         &60 & 20 & $ r_{20}=40, r_i=0 for i\ne20    $\\
$[21]$        &60 & 20 & $r_1=r_{20}=10, r_{10}=20 r_i=0$ for $ i\ne1,10,20$\\
&&&\\
$[22]$        &60 & 40 & $r_1=20, r_i=0  $ for $ i\ne1$\\
$[23]   $         &60 & 40 & $ r_{40}=20, r_i=0$ for $i\ne40   $\\
$[24]$        &60 & 40 & $r_{2i-1}=1,r_{2i}=0,$ for $ i=1,2,...,20$\\
&&&\\
$[25]$        &60 & 50 & $r_1=10, r_i=0  $ for $ i\ne1$\\
$[26]   $         &60 & 50 & $ r_{50}=10, r_i=0$ for $ i\ne1    $\\
$[27]$        &60 & 50 & $r_1=r_{50}=5, r_i=0 $ for $ i\ne1,50$\\

                      \hline
\end{tabular} 
\end{center}
\end{table}

\begin{table}
\begin{center}
\caption{ Simulated critical values of ${\rm H}_{m:n}$ }\label{t10}
\begin{tabular}{cccccc}
\hline
Scheme no. &  $ {\rm H}_{m:n}$  & Scheme no. & ${\rm H}_{m:n}$  & Scheme no. & ${\rm H}_{m:n}$\\
\hline
$[1]  $      &0.1069 & $[10]$ & 0.1033 & $[19]$ & 0.0633 \\
$[2]  $      &0.1062 & $[11]$ & 0.0941 & $[20]$ & 0.0595 \\
$[3]  $      &0.1060 & $[12]$ & 0.0971 & $[21]$ & 0.0621 \\
$[4]  $      &0.0802 & $[13]$ & 0.0588 & $[22]$ & 0.0351 \\
$[5]  $      &0.0793 & $[14]$ & 0.0573 & $[23]$ & 0.0358 \\
$[6]  $      &0.0846 & $[15]$ & 0.0602 & $[24]$ & 0.0370 \\
$[7]  $      &0.0646 & $[16]$ & 0.0424 & $[25]$ & 0.0296 \\
$[8]  $      &0.0661 & $[17]$ & 0.0431 & $[26]$ & 0.0300 \\
$[9]  $      &0.0671 & $[18]$ & 0.0425 & $[27]$ & 0.0298 \\

                      \hline
\end{tabular} 
\end{center}
\end{table}

In Table 8, we present the estimated power of the our proposed test, Balakrishnan et al.'s    \cite{bala2004}  T-statistic and  Pakyari  and Balakrishnan's    \cite{pak2013} test statistics when the null hypothesis stipulates normal and the alternative hypothesis corresponds to Student's t with  three and four degrees of freedom, Logistic distribution and double exponential distribution. From this table it is apparent that for a symmetric heavy-tailed alternative while testing for normality, the  test statistic, ${\rm H}_{m:n}$, that we have proposed, has possessed better power than Balakrishnan et al.'s    \cite{bala2004}  T-statistic and  Pakyari and Balakrishnan's    \cite{pak2013} 
test statistics in 78 out of 108 situations. Also, when $n=20$ in the Student's t distribution with three degrees of freedom in 4 out of 9 situations, in the Student's t distribution with four degrees of freedom in 5 out of 9, in the Logistic $(0,1)$ in 6 out of 9 situations and in the double exponential $(0,1)$ in 7 out of 9 situations, in the case $n=40$  in the Student's t distribution with three degrees of freedom in 7 out of 9 situations, in the Student's t distribution with four degrees of freedom in 7 out of 9, in the Logistic $(0,1)$ in 7 out of 9 situations and in the double exponential $(0,1)$ in 7 out of 9 situations  and in the cases $n=60$  in the Student's t distribution with three degrees of freedom in 6 out of 9 situations, in the Student's t distribution with four degrees of freedom in 8 out of 9, in the Logistic $(0,1)$ in 8 out of 9 situations and in the double exponential $(0,1)$ in 6 out of 9 situations our test statistic possessed better power than Balakrishnan et al.'s    \cite{bala2004}  T-statistic and  Pakyari  and Balakrishnan's    \cite{pak2013}  
test statistics. Also, in early censoring schemes ([1],[4],[7],[10],[13],[16],[19],[22],[25]), the $G_{m:n}^{(2)}$ statistic has the most power in 6 out 36 situations, the $G_{m:n}^{(3)}$ statistic has the most power in 11 out 36 situations, the $Q_{m:n}$ statistic has the most power in 1 out 36 situations, the $T_{m:n}^{(2)}$ statistic has the most power in 4 out 36 situations, the T-statistic has the most power in 4 out 36 situations and the ${\rm H}_{m:n}$ statistic has the most power in 10 out 36 situations. In addition, in non-early censoring schemes ([2],[3],[5],[6],[8],[9],[11],[12],[14],[15],[17],[18],[20],[21],[23],[24],[26],[27]), the $G_{m:n}^{(3)}$ statistic has the most power in 1 out 72 situations, the $Q_{m:n}$ statistic has the most power in 1 out 72 situations, the T-statistic has the most power in 2 out 72 situations and the ${\rm H}_{m:n}$ statistic has the most power in 68 out 72 situations. Also note that, as one would normally expect, it can be observed from the values in the Table 8 that the power increases as the degree of censoring $(1-m/n)$ decreases. Finally, based on this results and comprative findings, we recomend the use of $G_{m:n}^{(3)}$ and ${\rm H}_{m:n}$ statistics for the case of early censoring and the use of ${\rm H}_{m:n}$ statistic for the case of non-early censoring.

\begin{table}
\begin{center}
\begin{footnotesize}
{Table 8: Estimated powers for Student's t distribution with three and four degrees of freedom, Logistic  and double exponential distribution with parameters $\mu =0$ and $\sigma=1$  with $\alpha=0.1$. }\label{t8}
\begin{tabular}
{@{\hspace{1mm}}c@{\hspace{2mm}}@{\hspace{2mm}}c@{\hspace{2mm}}@{\hspace{2mm}}c@{\hspace{2mm}}c@{\hspace{2mm}}c@{\hspace{1mm}}c@{\hspace{2mm}}c@{\hspace{2mm}}c@{\hspace{2mm}}c@{\hspace{2mm}}c@{\hspace{1mm}}c@{\hspace{1mm}}c@{\hspace{1mm}}c@{\hspace{1mm}}c@{\hspace{1mm}}}
\hline
Scheme &  dist. &$C^+_{m:n}$  &$C^-_{m:n}$&$C_{m:n}$&$K_{m:n}$&$T^{(1)}_{m:n}$& $T^{(2)}_{m:n}$& $G_{m:n}$ & $Q_{m:n}$&$G^{(2)}_{m:n}$&$G^{(3)}_{m:n}$&$T$&${\rm H}_{m:n}$\\
no.   &&&&&&&&&&&&&\\
\hline
$[1]$&$t_{(3)}$	&0.2224&0.1112&0.2251&0.1681&0.2266&0.2158 &0.2438&0.2528&\textbf{0.2676}&0.2341&0.2567&0.2638\\
     &$t_{(4)}$	&0.1807&0.1016&0.1770&0.1303&0.1738&0.1661 &0.1946&0.2043&\textbf{0.2211}&0.1920&0.2056&0.2172\\
     & $L(0,1)$	  &0.1232&0.1003&0.1114&0.0944&0.1199&0.1122 &0.1369&0.1376&\textbf{0.1509}&0.1305&0.1217&0.1408\\
     &$DE(0,1)$	  &0.1912&0.0943&0.1803&0.1251&0.1725&0.1668 &0.2188&0.2322&\textbf{0.2377}&0.2216&0.1909&0.2286\\
     &&&&&&&&&&&&&\\
$[2]$&$t_{(3)}$	&0.2450&0.0507&0.1463&0.1459&0.1541&0.1554 &0.2057&0.2272&0.1579&0.0754&0.2461&\textbf{0.2895}\\
     &$t_{(4)}$	&0.1988&0.0616&0.1221&0.1234&0.1254&0.1257 &0.1678&0.1839&0.1335&0.0766&0.1965&\textbf{0.2380}\\
     &$L(0,1)$	&0.1440&0.0682&0.0937&0.0985&0.1069&0.0989 &0.1330&0.1321&0.1039&0.0956&0.1330&\textbf{0.1585}\\
     &$DE(0,1)$	&0.2222&0.0489&0.1198&0.1313&0.1310&0.1317 &0.1827&0.2078&0.1236&0.0694&0.1746&\textbf{0.2555}\\
     &&&&&&&&&&&&&\\
$[3]$&$t_{(3)}$	&0.2592&0.0471&0.1887&0.1288&0.1802&0.1771 &0.2316&0.2541&0.2883&0.2570&0.2543&\textbf{0.3000}\\
     &$t_{(4)}$	&0.2066&0.0572&0.1480&0.1098&0.1433&0.1385 &0.1871&0.2004&0.2276&0.2085&0.2018&\textbf{0.2382}\\
     &$L(0,1)$	&0.1458&0.0715&0.1094&0.0920&0.1157&0.1092 &0.1279&0.1345&0.1577&0.1364&0.1297&\textbf{0.1654}\\
     &$DE(0,1)$	&0.2438&0.0402&0.1671&0.1210&0.1774&0.1673 &0.2116&0.2423&0.2569&0.2372&0.1833&\textbf{0.2902}\\
     &&&&&&&&&&&&&\\
$[4]$&$t_{(3)}$	&0.2233&0.1464&0.2510&0.2000&0.2598&0.2583 &0.2523&0.2716&0.3036&0.2802&\textbf{0.3071}&0.3008\\
     &$t_{(4)}$	&0.1772&0.1200&0.1879&0.1505&0.1921&0.1916 &0.2021&0.2161&0.2424&0.2234&\textbf{0.2490}&0.2364\\
     &$L(0,1)$	&0.1235&0.0972&0.1318&0.1017&0.1180&0.1199 &0.1416&0.1471&\textbf{0.1748}&0.1669&0.1411&0.1673 \\
     &$DE(0,1)$	&0.1915&0.1195&0.2149&0.1647&0.2148&0.2029 &0.2332&0.2661&\textbf{0.2904}&0.2734&0.2340&0.2702\\
     &&&&&&&&&&&&&\\
$[5]$&$t_{(3)}$	&0.2630&0.0448&0.1848&0.1526&0.1892&0.1873 &0.2579&0.2807&0.3003&0.2657&0.3190&\textbf{0.3526}\\
     &$t_{(4)}$	&0.2036&0.0546&0.1442&0.1229&0.1462&0.1438 &0.2038&0.2227&0.2335&0.2112&0.2484&\textbf{0.2797}\\
     &$L(0,1)$	&0.1463&0.0721&0.1083&0.0929&0.1066&0.1136 &0.1393&0.1591&0.1585&0.1380&0.1311&\textbf{0.1927}\\
     &$DE(0,1)$	&0.2783&0.0343&0.1835&0.1485&0.1819&0.1890 &0.2483&0.3087&0.3112&0.2776&0.1980&\textbf{0.3682}\\
     &&&&&&&&&&&&&\\
$[6]$&$t_{(3)}$	&0.0952&0.2467&0.1784&0.1908&0.2123&0.2259 &0.2703&0.2895&0.2732&0.2211&\textbf{0.3386}&0.3367\\
     &$t_{(4)}$	&0.0819&0.1906&0.1373&0.1491&0.1626&0.1693 &0.2159&0.2277&0.2173&0.1803&0.2727&\textbf{0.2748}\\
     &$L(0,1)$	&0.0767&0.1309&0.1092&0.1013&0.1064&0.1076 &0.1497&0.1541&0.1435&0.1391&0.1434&\textbf{0.1826}\\
     &$DE(0,1)$	&0.0843&0.2266&0.1560&0.1698&0.1807&0.1941 &0.2406&0.2847&0.2605&0.2406&0.2051&\textbf{0.3234}\\
     &&&&&&&&&&&&&\\
$[7]$&$t_{(3)}$	&0.2183&0.1881&0.2718&0.2332&0.2894&0.2865 &0.2694&0.2832&0.2942&0.2586&\textbf{0.3407}&0.3193\\
     &$t_{(4)}$	&0.1698&0.1433&0.2018&0.1629&0.2102&0.2080 &0.2119&0.2216&0.2364&0.2076&\textbf{0.2755}&0.2645\\
     &$L(0,1)$	&0.1205&0.1086&0.1183&0.0984&0.1231&0.1249 &0.1484&0.1557&0.1540&0.1518&0.1539&\textbf{0.1707}\\
     &$DE(0,1)$	&0.1926&0.1646&0.2253&0.1824&0.2353&0.2485 &0.2507&0.2738&0.2740&0.2617&0.2668&\textbf{0.2850}\\
\end{tabular} 
\end{footnotesize}
\end{center}
\end{table}

\begin{table}[http]
\begin{center}
 \begin{footnotesize}
{Table 8: Continued }
\begin{tabular}
{@{\hspace{1mm}}c@{\hspace{2mm}}@{\hspace{2mm}}c@{\hspace{2mm}}@{\hspace{2mm}}c@{\hspace{2mm}}c@{\hspace{2mm}}c@{\hspace{1mm}}c@{\hspace{2mm}}c@{\hspace{2mm}}c@{\hspace{2mm}}c@{\hspace{2mm}}c@{\hspace{1mm}}c@{\hspace{1mm}}c@{\hspace{1mm}}c@{\hspace{1mm}}c@{\hspace{1mm}}}
\hline
Scheme &  dist. &$C^+_{m:n}$  &$C^-_{m:n}$&$C_{m:n}$&$K_{m:n}$&$T^{(1)}_{m:n}$& $T^{(2)}_{m:n}$& $G_{m:n}$ & $Q_{m:n}$&$G^{(2)}_{m:n}$&$G^{(3)}_{m:n}$&$T$&${\rm H}_{m:n}$\\
no.   &&&&&&&&&&&&&\\
\hline
$[8]$&$t_{(3)}$	&0.2533&0.0740&0.2055&0.1652&0.2048&0.2021 &0.2799&0.3089&0.3164&0.2917&0.3309&\textbf{0.3684}\\
     &$t_{(4)}$	&0.1963&0.0701&0.1581&0.1251&0.1564&0.1558 &0.2216&0.2406&0.2464&0.2266&0.2607&\textbf{0.2987}\\
     &$L(0,1)$	&0.1374&0.0761&0.1137&0.0954&0.1071&0.1090 &0.1505&0.1650&0.1668&0.1400&0.1391&\textbf{0.1900}\\
     &$DE(0,1)$	&0.2665&0.0625&0.2117&0.1632&0.2106&0.2097 &0.2842&0.3321&0.3240&0.3061&0.2100&\textbf{0.3941}\\
     &&&&&&&&&&&&\\
$[9]$&$t_{(3)}$	&0.1531&0.2330&0.2292&0.2185&0.2621&0.2663 &0.2839&0.2997&0.2805&0.2255&\textbf{0.3479}&0.3204\\
     &$t_{(4)}$	&0.1178&0.1722&0.1679&0.1584&0.1911&0.2226 &0.2381&0.2216&0.1805&\textbf{0.2810}&0.2607&0.2681\\
     &$L(0,1)$	&0.0944&0.1156&0.1169&0.1058&0.1124&0.1225 &0.1594&\textbf{0.1632}&0.1523&0.1300&0.1447&0.1620\\
     &$DE(0,1)$	&0.1326&0.1940&0.2005&0.1843&0.2153&0.2279 &0.2756&0.2964&0.2940&0.2189&0.2237&\textbf{0.2992}\\
     &&&&&&&&&&&&\\
$[10]$&$t_{(3)}$&0.3111&0.1037&0.3002&0.2300&0.3094&0.2918 &0.3258&0.3457&0.3408&\textbf{0.3507}&0.3098&0.3380\\
     &$t_{(4)}$	&0.2386&0.0880&0.2204&0.1637&0.2253&0.2111 &0.2521&0.2640&0.2750&\textbf{0.2759}&0.2432&0.2577\\
     &$L(0,1)$	&0.1514&0.0746&0.1295&0.1096&0.1328&0.1269 &0.1585&0.1581&0.1691&\textbf{0.1960}&0.1477&0.1584\\
     &$DE(0,1)$	&0.2549&0.0782&0.2242&0.1646&0.2329&0.2091 &0.2691&0.2883&0.2814&\textbf{0.3279}&0.2313&0.2556\\
 &&&&&&&&&&&&\\
$[11]$&$t_{(3)}$&0.3057&0.0404&0.1832&0.1939&0.1881&0.1952 &0.2153&0.2417&0.0773&0.0473&0.03170&\textbf{0.3527}\\
     &$t_{(4)}$	&0.2384&0.0488&0.1411&0.1539&0.1462&0.1520 &0.1728&0.1920&0.0753&0.0589&0.2646&\textbf{0.2847}\\
     &$L(0,1)$	&0.1541&0.0704&0.1024&0.1173&0.0967&0.1101 &0.1180&0.1301&0.0745&0.0855&0.1312&\textbf{0.1800}\\
     &$DE(0,1)$	&0.2041&0.0484&0.1156&0.14434&0.1178&0.1268 &0.1489&0.1594&0.0640&0.0540&0.1545&\textbf{0.2384}\\
     &&&&&&&&&&&&&\\   
$[12]$&$t_{(3)}$	&0.1814&0.0720&0.1291&0.1465&0.1715&0.1835 &0.2912&0.3369&0.3318&0.2963&0.3703&\textbf{0.4090}\\
     &$t_{(4)}$	&0.1413&0.0594&0.1048&0.1082&0.1263&0.1393 &0.2261&0.2566&0.2577&0.2303&0.2819&\textbf{0.3171}\\
     &$L(0,1)$	&0.1091&0.0614&0.0841&0.0842&0.897&0.0875 &0.1341&0.1686&0.1528&0.1483&0.1453&\textbf{0.2013}\\
     &$DE(0,1)$	&0.1527&0.0372&0.0987&0.0914&0.1201&0.1205 &0.2096&0.2754&0.2429&0.2413&0.2469&\textbf{0.3116}\\
     &&&&&&&&&&&&&\\
$[13]$&$t_{(3)}$	&0.1363&0.2004&0.3627&0.3282&0.3897&0.3860 &0.3567&0.3866&0.4212&\textbf{0.4415}&0.3797&0.4318\\
     &$t_{(4)}$	&0.2544&0.1407&0.2617&0.2208&0.2738&0.2725 &0.2689&0.2940&0.3319&\textbf{0.3501}&0.3018&0.3369\\
     &$L(0,1)$	&0.1415&0.1015&0.1456&0.1200&0.1407&0.1471 &0.1645&0.1700&0.1828&\textbf{0.2236}&0.1651&0.2040\\
     &$DE(0,1)$	&0.2980&0.1557&0.2974&0.2668&0.3148&0.3241 &0.2969&0.33276&0.3497&\textbf{0.3984}&0.2635&0.3751\\
     &&&&&&&&&&&&&\\
$[14]$&$t_{(3)}$	&0.4075&0.0375&0.2934&0.2543&0.3087&0.3124 &0.3095&0.3592&0.3960&0.3913&0.4820&\textbf{0.5050}\\
     &$t_{(4)}$	&0.3070&0.0395&0.2076&0.1807&0.2100&0.2145 &0.2359&0.2706&0.3071&0.3049&0.3729&\textbf{0.3886}\\
     &$L(0,1)$	&0.1848&0.0623&0.1267&0.1070&0.1177&0.1229 &0.1504&0.1612&0.1758&0.1968&0.1374&\textbf{0.2325}\\
     &$DE(0,1)$	&0.4194&0.0186&0.2619&0.2113&0.2737&0.2944 &0.2724&0.3205&0.3443&0.3652&0.2506&\textbf{0.4484}\\
     &&&&&&&&&&&&&\\
$[15]$&$t_{(3)}$&0.1300&0.3242&0.2457&0.2549&0.2921&0.3088&0.3628 &0.4177&0.4143&0.4054&0.4822&\textbf{0.5062}\\
     &$t_{(4)}$	&0.0962&0.2350&0.1700&0.1765&0.2003&0.2096 &0.2773&0.3184&0.3150&0.3109&0.3797&\textbf{0.3987}\\
     &$L(0,1)$	&0.0751&0.1377&0.1073&0.1152&0.1136&0.1211 &0.1723&0.1905&0.1932&0.1932&0.1624&\textbf{0.2454}\\
     &$DE(0,1)$	&0.1071&0.2996&0.2098&0.2307&0.2793&0.2967 &0.3308&0.4087&0.3884&0.4057&0.2597&\textbf{0.4997}\\
     &&&&&&&&&&&&&\\
\end{tabular} 
 \end{footnotesize}
\end{center}
\end{table}

\begin{table}[htp]
\begin{center}
     \begin{footnotesize}
{Table 8: Continued }
\begin{tabular}
{@{\hspace{1mm}}c@{\hspace{2mm}}@{\hspace{2mm}}c@{\hspace{2mm}}@{\hspace{2mm}}c@{\hspace{2mm}}c@{\hspace{2mm}}c@{\hspace{1mm}}c@{\hspace{2mm}}c@{\hspace{2mm}}c@{\hspace{2mm}}c@{\hspace{2mm}}c@{\hspace{1mm}}c@{\hspace{1mm}}c@{\hspace{1mm}}c@{\hspace{1mm}}c@{\hspace{1mm}}}
\hline
Scheme &  dist. &$C^+_{m:n}$  &$C^-_{m:n}$&$C_{m:n}$&$K_{m:n}$&$T^{(1)}_{m:n}$& $T^{(2)}_{m:n}$& $G_{m:n}$ & $Q_{m:n}$&$G^{(2)}_{m:n}$&$G^{(3)}_{m:n}$&$T$&${\rm H}_{m:n}$\\
no.   &&&&&&&&&&&&&\\
\hline
$[16]$&$t_{(3)}$	&0.3546&0.2936&0.4107&0.3990&0.4461&0.4501 &0.3786&0.4153&0.4316&0.4186&0.4152&\textbf{0.4737}\\
     &$t_{(4)}$	&0.2499&0.2000&0.2784&0.2635&0.3096&0.3087 &0.2802&0.3084&0.3285&0.3197&0.3281&\textbf{0.3666}\\
     &$L(0,1)$	&0.1411&0.1157&0.1412&0.1272&0.1521&0.1467 &0.1626&0.1852&0.1987&0.2013&0.1903&\textbf{0.2168}\\
     &$DE(0,1)$	&0.3175&0.2432&0.3443&0.3455&0.3900&0.3992 &0.3292&0.3717&0.3921&0.4011&0.2906&\textbf{0.4573}\\
     &&&&&&&&&&&&&\\
$[17]$&$t_{(3)}$	&0.3961&0.0941&0.3201&0.2794&0.3346&0.3410 &0.3847&0.4430&0.4615&0.4671&0.4998&\textbf{0.5470}\\
     &$t_{(4)}$	&0.2920&0.0661&0.2259&0.1878&0.2373&0.2387 &0.2865&0.3333&0.3481&0.3508&0.3894&\textbf{0.4317}\\
     &$L(0,1)$	&0.1695&0.0586&0.1278&0.1057&0.1246&0.1386 &0.1648&0.1950&0.2115&0.2326&0.1632&\textbf{0.2567}\\
     &$DE(0,1)$	&0.4396&0.0527&0.3429&0.2903&0.3656&0.3753 &0.3772&0.4531&0.4757&0.5266&0.2870&\textbf{0.5688}\\
     &&&&&&&&&&&&&\\
$[18]$&$t_{(3)}$	&0.3741&0.1213&0.3201&0.2835&0.3379&0.3416 &0.3863&0.4400&0.4490&0.4538&0.4439&\textbf{0.5330}\\
     &$t_{(4)}$	&0.2769&0.0849&0.2244&0.1925&0.2365&0.2402 &0.2889&0.3318&0.3426&0.3452&0.3461&\textbf{0.4115}\\
     &$L(0,1)$	&0.1582&0.0726&0.1246&0.1061&0.1395&0.1321 &0.1864&0.1985&0.2037&0.2036&0.1576&\textbf{0.2479}\\
     &$DE(0,1)$	&0.3883&0.0976&0.3167&0.2919&0.3603&0.3329 &0.3677&0.4269&0.4321&0.4699&0.2503&\textbf{0.5099}\\
     &&&&&&&&&&&&&\\
$[19]$&$t_{(3)}$	&0.4230&0.1869&0.4158&0.3865&0.4521&0.4469 &0.4163&0.4513&0.4404&\textbf{0.5027}&0.3904&0.4785\\
     &$t_{(4)}$	&0.3098&0.1264&0.2962&0.2576&0.3154&0.3087 &0.3096&0.3428&0.3365&\textbf{0.3904}&0.3073&0.3710\\
     &$L(0,1)$	&0.1777&0.0845&0.1550&0.1286&0.1536&0.1423 &0.1826&0.1988&0.2162&\textbf{0.2221}&0.1680&0.2102\\
     &$DE(0,1)$	&0.3499&0.1380&0.3285&0.2934&0.3425&0.3261 &0.3390&\textbf{0.3982}&0.3413&0.3860&0.2666&0.3971\\
     &&&&&&&&&&&&&\\
$[20]$&$t_{(3)}$	&0.4547&0.0220&0.3208&0.2860&0.3392&0.3472 &0.2853&0.3377&0.2925&0.2383&0.5208&\textbf{0.5308}\\
     &$t_{(4)}$	&0.3478&0.0312&0.2229&0.2003&0.2332&0.2406 &0.2152&0.2517&0.2161&0.1713&0.4021&\textbf{0.4154}\\
     &$L(0,1)$	&0.1985&0.0512&0.1223&0.1136&0.1235&0.1286 &0.1461&0.1599&0.1195&0.1085&0.1399&\textbf{0.2434}\\
     &$DE(0,1)$	&0.3416&0.0224&0.2026&0.1761&0.2149&0.2176 &0.2058&0.2468&0.1533&0.1207&0.2026&\textbf{0.3815}\\
          &&&&&&&&&&&&&\\
$[21]$&$t_{(3)}$	&0.1537&0.2744&0.2065&0.2770&0.2871&0.3231 &0.3994&0.4680&0.4739&0.4303&0.5826&\textbf{0.6001}\\
     &$t_{(4)}$	&0.1051&0.1762&0.1354&0.1757&0.1839&0.2116 &0.2988&0.3540&0.3527&0.3163&0.4496&\textbf{0.4700}\\
     &$L(0,1)$	&0.0830&0.0999&0.0826&0.0891&0.1023&0.1047 &0.1772&0.2033&0.1978&0.1796&0.1595&\textbf{0.2718}\\
     &$DE(0,1)$	&0.1355&0.1732&0.1407&0.1974&0.2154&0.2346  &0.3432&0.4297&0.4373&0.4069&0.2931&\textbf{0.5277}\\
     &&&&&&&&&&&&&\\
\end{tabular} 
     \end{footnotesize}
\end{center}
\end{table}

\begin{table}[hp]
\begin{center}
   \begin{footnotesize}
{Table 8: Continued }
\begin{tabular}
{@{\hspace{6mm}}c@{\hspace{2mm}}@{\hspace{2mm}}c@{\hspace{2mm}}@{\hspace{2mm}}c@{\hspace{2mm}}c@{\hspace{2mm}}c@{\hspace{1mm}}c@{\hspace{2mm}}c@{\hspace{2mm}}c@{\hspace{2mm}}c@{\hspace{2mm}}c@{\hspace{1mm}}c@{\hspace{1mm}}c@{\hspace{1mm}}c@{\hspace{1mm}}c@{\hspace{1mm}}}
\hline
Scheme &  Dis. &$C^+_{m:n}$  &$C^-_{m:n}$&$C_{m:n}$&$K_{m:n}$&$T^{(1)}_{m:n}$& $T^{(2)}_{m:n}$& $G_{m:n}$ & $Q_{m:n}$&$G^{(2)}_{m:n}$&$G^{(3)}_{m:n}$&$T$&${\rm H}_{m:n}$\\
no.   &&&&&&&&&&&&\\
\hline
$[22]$&$t_{(3)}$	&0.4571&0.3707&0.5185&0.5201&0.5653&\textbf{0.5719}&0.4529&0.5039&0.5243&0.5370&0.4446&0.5715\\
     &$t_{(4)}$	&0.1963&0.0701&0.1581&0.1251&0.1564&0.1558 &0.2216&0.2406&0.2464&0.2266&0.2607&\textbf{0.4462}\\
     &$L(0,1)$	&0.1503&0.1237&0.1621&0.1491&0.1786&0.1764 &0.1805&0.2000&0.2090&0.2298&0.1911&\textbf{0.2441}\\
     &$DE(0,1)$	&0.3981&0.2911&0.4360&0.4678&0.5049&\textbf{0.5090} &0.3715&0.4322&0.4402&0.4912&0.3024&0.4984\\
     &&&&&&&&&&&&&\\
$[23]$&$t_{(3)}$	&0.5234&0.1172&0.4238&0.3831&0.4550&0.4567 &0.4385&0.5148&0.5779&0.5981&0.6351&\textbf{0.6714}\\
     &$t_{(4)}$	&0.3841&0.0652&0.2840&0.2454&0.3064&0.3083 &0.3228&0.3847&0.4409&0.4929&0.2607&\textbf{0.5298}\\
     &$L(0,1)$	&0.2163&0.0518&0.1472&0.1265&0.1473&0.1502 &0.1886&0.2124&0.2422&0.2622&0.1712&\textbf{0.3148}\\
     &$DE(0,1)$	&0.6119&0.0461&0.4696&0.4136&0.5146&0.5225 &0.4309&0.5162&0.6080&0.6556&0.4060&\textbf{0.6967}\\
     &&&&&&&&&&&&&\\
$[24]$&$t_{(3)}$	&0.3141&0.4451&0.4329&0.4479&0.4833&0.4927 &0.4523&0.5161&0.5161&0.5193&0.5263&\textbf{0.6261}\\
     &$t_{(4)}$	&0.1950&0.3175&0.2862&0.2975&0.3233&0.3324 &0.3371&0.3890&0.3876&0.3912&0.4167&\textbf{0.4838}\\
     &$L(0,1)$	&0.0990&0.1656&0.1488&0.1408&0.1557&0.1548 &0.1984&0.2225&0.2316&0.2301&0.1749&\textbf{0.2938}\\
     &$DE(0,1)$	&0.2946&0.4103&0.4173&0.4244&0.4865&0.4846 &0.4079&0.4856&0.4971&0.5211&0.3063&\textbf{0.5832}\\
     &&&&&&&&&&&&&\\
$[25]$&$t_{(3)}$	&0.4838&0.4407&0.5496&0.5684&0.6065&\textbf{0.6245}&0.4657&0.5189&0.5224&0.5280&0.4634&0.5940\\
     &$t_{(4)}$	&0.3296&0.2879&0.3724&0.3788&0.4189&0.4261 &0.3361&0.3809&0.3902&0.3910&0.36369&\textbf{0.4659}\\
     &$L(0,1)$	&0.1552&0.1360&0.1606&0.1605&0.1834&0.1788 &0.1849&0.1992&0.2157&0.2224&0.1977&\textbf{0.2639}\\
     &$DE(0,1)$	&0.4304&0.3850&0.4865&0.5501&0.5685&\textbf{0.5850} &0.3942&0.4436&0.4715&0.4904&0.3039&0.5128\\
     &&&&&&&&&&&&&\\
$[26]$&$t_{(3)}$	&0.4908&0.1752&0.4140&0.4005&0.4492&0.4552 &0.4645&0.5355&0.5497&0.5744&0.5754&\textbf{0.6567}\\
     &$t_{(4)}$	&0.3548&0.1029&0.2788&0.2605&0.3021&0.3077 &0.3456&0.3989&0.4124&0.4333&0.4447&\textbf{0.5107}\\
     &$L(0,1)$	&0.1839&0.0697&0.1562&0.1206&0.1558&0.1527 &0.1997&0.2219&0.2333&0.2504&0.1666&\textbf{0.3074}\\
     &$DE(0,1)$	&0.5219&0.1568&0.4641&0.4382&0.5160&0.4952 &0.4495&0.5341&0.5499&0.6007&0.3264&\textbf{0.6580}\\
     &&&&&&&&&&&&&\\
$[27]$&$t_{(3)}$	&0.4670&0.2110&0.4169&0.4109&0.4574&0.4616 &0.4619&0.5236&0.5275&0.5526&0.4998&\textbf{0.6217}\\
     &$t_{(4)}$	&0.3330&0.1286&0.2819&0.2719&0.3149&0.3146 &0.3428&0.3917&0.3980&0.4155&0.3871&\textbf{0.4983}\\
     &$L(0,1)$	&0.1796&0.0861&0.1443&0.1315&0.1537&0.1545 &0.1905&0.2172&0.2170&0.2393&0.1635&0.\textbf{2842}\\
     &$DE(0,1)$	&0.4800&0.2140&0.4486&0.4521&0.5073&0.5037 &0.4226&0.4995&0.5100&0.5325&0.2744&\textbf{0.6037}\\
     &&&&&&&&&&&&&\\
                      \hline
\end{tabular} 
\end{footnotesize}
\end{center}
\end{table}

\section{Illustrative data analyses}
In this section, the wire connection strength data from Nelson  \cite{nelson}, (Table 5.1, p. 111) are considered. these data, originally studied by King  \cite{king}, concern the breaking strength of 23
wire connections. The wires were bonded at one end to a semiconductor wafer and at the other
end to a terminal post. The first two and the last one of the observations were eliminated from
the analysis due to validity suspection of the data; see Nelson \cite{nelson}, for more details.
 Pakyari  and Balakrishnan    \cite{pak2013} randomly generated a progressively Type-II censored sample of size $m = 10$ from $n = 20$ observations.
Table 9 presents the data and the corresponding progressive censoring scheme.
The possibility of fitting a normal model to the data was done by  Nelson  \cite{nelson}, and we, therefore, tested for normality.
Table 10 presents the test statistics and their corresponding p-values. The normal model is strongly supported by all the test statistics for describing the wire connection strength data.  Results in Table 9 show that the p-value of the our test statistic, ${\rm H}_{m:n}$, is  greater than other p-values.

\begin{table}[h]
\begin{center}
{Table 9: Wire connection strength data and the progressive Type-II censoring scheme. }\label{20}
\begin{tabular}{ccccccccccc}
\hline
$i$ & 1 & 2&3&4&5&6&7&8&9&10\\
$x_{i:m:n}$ &550&750&950&1150&1150&1150&1350&1450&1550&1850\\
$r_i$&0&2&1&0&3&0&0&2&0&2\\
                      \hline
\end{tabular} 
\end{center}
\end{table}
 
\begin{table}[h]
\begin{center}
{Table 10: Test statistics and the corresponding p-values for the data given
in Table 9 when testing for the normal distribution. }\label{t20}
\begin{tabular}{ccccccc}
\hline
Criterion  &$C^+_{m:n}$  &$C^-_{m:n}$&$C_{m:n}$&$K_{m:n}$&$T^{(1)}_{m:n}$& $T^{(2)}_{m:n}$\\
Test statistic &0.0946&0.0893&0.0946&0.1839&0.0021&0.0352\\
p-value&0.6576&0.3809&0.7057&0.5364&0.8020&0.8735\\
&&&&&&\\ 
 Criterion & $G_{m:n}$ & $Q_{m:n}$&$G^{(2)}_{m:n}$&$G^{(3)}_{m:n}$&$T$&${\rm H}_{m:n}$ \\
Test statistic &6.8499&10.9208&26.7465&63.8562&0.4568 &0.3220 \\
p-value &0.7152&0.6476&0.6879&0.6689&0.6450 &0.8091\\
                      \hline1
\end{tabular} 
\end{center}
\end{table}
 
\section{Conclusions}
In this paper, we proposed a simple and powerful test for normality based on progressively Type-II censored order statistics and compared this new test with all previous tests proposed for normality. Using a simulation study, consistency of our test was illustrated and also power of the test for some various alternatives were obtained and summarized.   
It was apparent from Table 6 that none of the tests considered performs better than all other tests against all alternatives.  Comparing  with  other tests, however, the proposed test ${\rm H}_{m:n}$, was the most powerful with respect to approximately all censoring schemes. Then, the performance of our test was examined for a real data set and the results were completely coincided with the other tests.





\bibliographystyle{elsarticle-num}
\bibliography{<your-bib-database>}



\end{document}